\input amstex 
\hsize=30truecc
\documentstyle{amsppt} 
\baselineskip=12pt 
\pagewidth{30pc}
\pageheight{52pc}
\loadbold
\TagsOnRight
\nologo
\nopagenumbers

\define\na{\Bbb N} 
\define\proc{$(\Cal P^\z,\mu,\sigma)$}

\define\procbr{$((\Cal P_B^r)^\z,\mu_B,\sigma_B)$}
 
\define\sq{sequence} 

\define\z{\Bbb Z}
\define\r{\Bbb R}

\topmatter  
 
\title The law of series\endtitle
\author 
T. Downarowicz\footnote{corresponding author: downar\@\,im.pwr.wroc.pl\hfill\break} and Y. Lacroix
\endauthor 
\address 
Institute of Mathematics and Computer Science, Wroclaw University of
Technology, Wy\-brze$\dot{\text z}$e Wys\-pia{\'n}\-skie\-go 27, 
50@-370 Wroc{\l}aw, 
Poland
\endaddress
\date May 16, 2008\enddate
\email 
downar\@pwr.wroc.pl
\endemail
\thanks
This paper was written during the first author's visits at CPT/ISITV in 2005 and 2006, 
supported by CNRS and ISITV. Research of the first author is supported from resources 
for science in years 2005-2008 as research project (grant MENII 1 P03A 021 29, Poland)  
\endthanks 
\address 
Institut des Sciences de l'Ing\'enieur de Toulon et du Var,
Avenue G. Pompidou, B.P. 56, 83162 La Valette du Var Cedex,
France
\endaddress
\email 
yves.lacroix\@univ-tln.fr
\endemail
\subjclass 37A50, 37A35, 37A05, 60G10
\endsubjclass
\keywords stationary random process, positive entropy, return time
statistic,  hitting time statistic, repelling, attracting, limit law, the
law of series, typical property 
\endkeywords

\abstract 
We consider an ergodic process on finitely many states, with positive entropy. Our first main result asserts that the distribution 
function of the normalized waiting time for the first visit to a small (i.e., over a long block) cylinder set $B$
is, for majority of such cylinders and up to epsilon, dominated by the exponential distribution function $1-e^{-t}$. That is,
the occurrences of so understood ``rare event'' $B$ along the time axis can appear either with gap sizes of nearly exponential distribution (like in the independent Bernoulli process), or they ``attract'' each-other. Our second main result states 
that a {\it typical} ergodic process of positive entropy has the following property: the distribution functions of the normalized 
hitting times for the majority of cylinders $B$ of lengths $n'$ converge to zero along a \sq\ $n'$ whose upper density is 1. 
The occurrences of such a cylinder $B$ ``strongly attract'', i.e., they appear in ``series'' of many frequent repetitions 
separated by huge gaps of nearly complete absence. 

These results, when properly and carefully interpreted, shed some new light, in purely statistical terms, independently 
from physics, on a century old (and so far rather avoided by serious science) common-sense phenomenon known as {\it the 
law of series}, asserting that rare events in reality, once occurred, have a mysterious tendency for untimely repetitions. 
\endabstract 

\endtopmatter 
\document

\comment
\heading Note \endheading
This paper resulted from studying asymptotic laws for return/hitting time statistics in stationary
processes, a field in ergodic theory rapidly developing in the recent years (see e.g. [A-G], [C], [C-K], 
[D-M], [H-L-V], [L] and the reference therein). Our first result significantly contributes to this area 
due to both its generality and strength. After having completely written the proof, during a free-minded discussion, 
the authors have discovered an astonishing interpretation of the result, clear even in terms of the common 
sense understanding of random processes. The organization of the paper is aimed to emphasise this discovery. 
The consequences for the field of asymptotic laws are moved toward the end of the paper. 
\endcomment

\heading Introduction\endheading 
We study the distribution functions of the hitting (and automatically also return) time statistics for small cylinder 
sets in processes on finitely symbols. We refer the reader to the rich literature on the subject (e.g. [A-G], [C], 
[C-K], [D-M], [H-L-V], [L] and the reference therein) for the recent developments in this field. Many works concentrate
on determining whether a process (or a class of processes) has ``exponential asymptotics'' or not. These attempts 
were successful in rather restricted classes of processes. Our Theorem 1 (and its variant, Theorem 3) is the first 
fully general result saying something concrete about all ergodic positive entropy processes, from this point of view. 
Namely, we prove that in such processes any essential limit distribution function for the hitting times is  majorized 
by the exponential law $1-e^{-t}$. In particular, this excludes many behaviors proved to exist in zero entropy, 
such as the presence of an essential limit law for the return times concentrated away from zero. 

This theorem sheds a new light on the extensively studied class of ergodic processes with positive entropy, where one 
could expect, all general properties have been established already long ago. It is impossible not to mention 
here the theorem of Ornstein and Weiss [O-W2] which relates the return times of long blocks to entropy. However, this 
theorem says nothing about the asymptotics of the distribution of the return times, because the logarithmic limit 
appearing in the statement is insensitive to the proportions between the gap sizes. 

Our approach is slightly different from the one represented in most papers on the return/hitting time assymptotics, 
as we are not interested in computing the limit laws ``at points'', i.e., along cylinders shrinking to a point $x$, 
where $x$ usually belongs to a positive (or full) measure set. We describe the restrictions on the distributions valid 
for ``majority'' of long cylinders $B$. The passage from our approach to the limit laws at points is described in the last section. 

The proof of Theorem 1 is rather complicated, yet entirely contained within the classics of ergodic theory; it relies on basic facts on entropy 
for partitions and sigma-fields, some elements of the Ornstein theory ($\epsilon$-independence), the Shannon-McMillan-Breiman 
Theorem, the Ornstein-Weiss Theorem on return times, the Ergodic Theorem, basics of probability and calculus. 

\smallskip
Our Theorem 2 belongs to the category describing typical (or generic) properties. It states that a typical ergodic process with positive 
entropy (see the last paragraph of this section for the meaning of typicality among positive entropy processes) has the following property 
which we call {\it strong attracting}\,: there exists a subsequence of lengths $(n')$ of upper density 1 in $\na$, such that the distribution 
functions of the normalized hitting times for the majority of cylinders $B$ of lengths $n'$ are ``flat'', i.e., close to zero on a long 
interval. Recall that only not long ago ([C-K]) it was 
discovered that some mixing (but still of entropy zero) transformations admit nonexponential asymptotics. Our result shows 
that even some Bernoulli processes do so, which, in particular, answers in the negative a question of Zaqueu Coelho [C]. 
\smallskip
Both inequalities between the distribution function of the normalized hitting time for an event $B$ and the exponential law $1-e^{-t}$ 
have nice and clear interpretations in terms of what we call {\it attracting} -- the tendency of the occurrences of $B$ to appear in series, and {\it repelling} -- the opposite tendency, toward a more uniform distribution of occurrences along the time axis. 
To our knowledge, these interpretations have not been addressed or discussed in any papers in the field. In these terms, our results 
can be expressed as follows: Theorem 1 -- in any positive entropy process the repelling of almost every sufficiently long cylinder 
$B$ is at most marginal; Theorem~2 -- within any measure-preserving system of positive entropy, if we ``draw'' a finite partition, then most likely it will generate a process, where nearly all long blocks of certain lengths (belonging to a large subset of $\na$) strongly attract. 

If we extrapolate this to processes and rare events running in reality, we obtain an astonishing contribution to the century 
old discussion about the so-called {\it law of series} (see the next section for more details).
\smallskip

Our understanding of typicality is somewhat different from the often considered setup, in which the set of all measure-preserving 
transformations (the automorphism group) on a fixed probability space is endowed with the topology of the weak convergence. In 
this setup, a typical transformation has entropy zero ([Ro]). Besides, the property we want to examine (strong attracting) depends on the generating partition, so we need to allow the partition to vary. Thus, we fix a measure-preserving system of positive entropy and $m\ge 2$, we consider all factor-processes generated by varying partitions into at most $m$ elements, and we adopt the notion of typicality with respect to the usual Rokhlin metric for partitions (which is complete on such partitions). Here, a typical process has positive entropy, so this approach is reasonable for studying ``typical properties of positive entropy systems''. 
Although we define typicality within a fixed system, strong attracting turns out to be typical inside every positive entropy 
system, which makes our notion of typicality for this property universal. 
\medskip
\centerline{ \it Acknowledgments }
\smallskip
The authors would like to thank Dan Rudolph for a hint leading to the first example of a positive entropy process with attracting,
and, in effect, to the discovery of the attracting/repelling asymmetry. (The same basic idea is still used in the proof of Theorem 
2.) We are grateful to Jean-Paul Thouvenot for his interest in the subject, substantial help, and the challenge to find a purely 
combinatorial proof of Theorem~1 (which we save for the future). 


\heading The common sense {\eightpoint \it LAW OF SERIES} versus our results \endheading

A ``series'' is noted in the every-day life, when a random event considered extremely rare 
happens more than once in a relatively short period of time. 
In the common sense, the {\it law of series} asserts that such series occur more often than they intuitively should,
indicating the existence of an unexplained physical force or statistical rule provoking them. For example, runs of 
good luck happen to gamblers, leading to high winnings (see [Wi] for the famous case of Charles Wells), people 
experience repetitions of similar unlucky events (hence the proverb``misfortune never comes alone''), or notice
series of strange coincidences without particular consequence, such as meeting people with the same last name on
the same day, seeing several times the same combination of digits in unrelated situations, etc.

An Austrian biologist dr. Paul Kammerer (1880-1926) was the first scientist to study this law. Although his book 
[Km] has attracted a lot of attention with its numerous suggestive examples, the scientific value of his 
``statistical'' interpretation is rather questionable. Kammerer himself lost authority due to accusations 
of manipulating his (unrelated to our topic) biological experiments.

Also some very serious scientists such as Swiss professor of philosophy Karl Gustav Jung (1875-1961), and a Nobel prize winner 
in physics, Austrian, Wolfgang Pauli (1900-1958), fascinated by examples of ``meaningful coincidences'' conjectured the existence
of undiscovered and mysterious ``attracting'' forces driving objects that are alike, or have common features, closer together in 
time and space, for which they coined a term ``synchronicity''. This includes attracting of repetitions of rare events in time, 
i.e., the {\it law of series}. Critics of synchronicity claim that all such ``unbelievable coincidencies'' and ``series'' occur 
at the rate complying with the statistics of pure randomness (see e.g. [Mi]). Human memory is keen to register them as more 
frequent simply because they are more distinctive. 

To be precise, let us agree that an event repeats in time by ``pure chance'' when it follows a Poisson process. 
In a typical realization of such a process, the distribution of signals along the time axis reveals a natural 
tendency to create spontaneous clusters, which can be easily taken for series, but are in fact just a feature of the random 
(unbiased) behavior. In order to say that some signal process obeys the {\it law of series}, one should detect in this 
process a tendency to create clusters stronger than in the Poisson process. It is possible to formally define such tendency 
without referring to the multidimensional distributions of the process, only to the 
single distribution of the normalized waiting time $V$ for the first signal. Because the waiting time for a signal in 
the Poisson process is exponential, such definition reduces to a simple inequality between the distribution function 
of $V$ and the function $1-e^{-t}$. This is exactly how we define attracting (see the next section). The extreme form of 
attracting, strong attracting, as we will define it, takes place when the signals occur in long series of frequency 
much higher than the probability of the signal, compensated by much longer periods of nearly complete absence.

$$
\align 
\text{\eightpoint repelling \ }& ..{...B}......\underline{B......B...B.....B}...{...B}.....
\underline{B....B}.....\underline{B....B..B}....B...{...B}..\\
\text{\eightpoint unbiased \ }& 
..{...B}........\underline{B....B..B....B}.....{...B}......\underline{B..B}
.......\underline{B..B.B}...{...B}...{...B}..\\
\text{\eightpoint attracting \ }&..{...B}..........\underline{B..B.B..B}........
{...B}.......\underline{BB}.........\underline{B.BB}....{...B}...{...B}..\\
\vspace {5pt}
\text{\eightpoint strong attr. \ }
&......B.BBBB.....................................................BB.BBB..B.B..\\
\endalign
$$ 
{\eightpoint\it Figure 1: Comparison between unbiased, repelling, attracting and strongly attracting distributions of occurrences 
of an event $B$ along the time.}
\bigskip\noindent

In many processes in reality, attracting or even strong attracting is perfectly understandable as a result of physical 
dependence. For example, many events reveal increased frequency of occurrences in so-called periods of propitious conditions, 
which in turn, follow a slowly changing steering process (e.g., floods following the climate changes). Such attracting, 
of course, is not the subject of the mystery behind the {\it law of series}. The challenge is to understand attracting for these 
events, for which we see no physical dependence and which are expected to have the unbiased behavior.

\medskip
With slight abuse of the complexity of life, our theorems can be interpreted to support the {\it law of series} as 
predominance of attracting for certain type of events. Reality is a realization of a huge measure-preserving system 
(obviously of positive entropy). Because we consider a single realization, we may assume ergodicity (a realization 
of a non-ergodic process belongs, almost surely, to an ergodic component). An ``elementary rare event'' whose occurrences 
cannot be fully predicted is a small cylinder set depending on a nondeterministic (i.e., also of positive entropy) 
factor-process generated by some finite partition of the phase space of this huge system. Then the majority of such 
elementary rare events reveal tendency to create series at least as strong as in the Poisson process (unbiased), or 
stronger. And in most cases this tendency will be in fact much stronger. Even if a real process is theoretically modeled 
by the Bernoulli process with an independent generator, so it is supposed to be unbiased (for example the process
of coin tosses), in reality the independent partition is always slightly perturbed, and then, by the typicality result, 
there will be an essential set of lengths $n'$ such that nearly all blocks of these lengths strongly attract. Because 
by Theorem 1, blocks of other lengths cannot essentially repel, ``in the average'', we will be dealing (against the 
intuition) with a substantial predominance of attracting for long configurations. 

Notice that the attracting is explained in purely statistical terms, without needing to understand the physical 
nature of the tiny dependencies in the perturbed generator. 

Of course, this hardly applies to gambling, because the event of, say, drawing a winning hand, is not a single cylinder, 
and it involves blocks probably too short for the attracting to take effect. But the theory may apply to some 
rare events in computer sciences, genetics or in other areas.

\heading Rigorous definitions and statements \endheading
We establish the notation necessary to formulate the main results. Let \proc\ be an ergodic process on 
finitely many symbols, i.e., $\#\Cal P<\infty$, $\sigma$ is the standard left shift map and $\mu$ is an 
ergodic shift-invariant probability measure on $\Cal P^\z$. Most of the time, we will identify finite blocks 
with their 
cylinder sets, i.e., we agree that $\Cal P^n=\bigvee_{i=0}^{n-1}\sigma^{-i}(\Cal P)$. Depending on the context, 
a block $B\in\Cal P^n$ is attached to some coordinates or it represents a ``word'' which may appear in different 
places along the $\Cal P$-names. We will also use the probabilistic language of random variables. Then 
$\mu\{R\in A\}$ ($A\subset\r$) will abbreviate $\mu(\{x\in\Cal P^\z:R(x)\in A\})$. Recall, that if the random 
variable $R$ is nonnegative and $F(t)=\mu\{R\le t\}$ is its distribution function, then the expected value 
of $R$ equals $\int_0^\infty 1-F(t)\,dt$.

For a set $B$ of positive measure let $R_B$ and $\overline R_B$ denote the random variables defined on $B$ 
(with the conditional measure $\mu_B =\frac\mu{\mu(B)}$) as the absolute and normalized first return time to $B$, 
respectively, i.e.,
$$
R_B(y) = \min\{i>0, \sigma^i(y)\in B\}, \ \ \overline R_B(y) = \mu(B)R_B(y).
$$
We denote by $\tilde F_B(t)$ the distribution function of $\overline R_B$.
Notice that, by the Kac Theorem ([Kc]), the expected value of $R_B$ equals $\frac1{\mu(B)}$, hence that of 
$\overline R_B$ is 1 (that is why we call it ``normalized''). We also define
$$
G_B(t) = \int_0^t 1-\tilde F_B(s)\, ds.
$$
Clearly, $G_B(t)\le \min\{t,1\}$ and the equality holds when $\tilde F_B(t) = 1_{[1,\infty)}$, that is, when $B$ 
occurs precisely with equal gaps (i.e., periodically); the gap size then equals $\frac1{\mu(B)}$.

Similarly, let $V_B$ be the random variable defined on $\Cal P^\z$ as the {\it hitting time statistic}, i.e., the waiting 
time for the first visit in $B$ (the defining formula is the same as for $R_B$, but this time it is regarded on the whole 
space with the measure $\mu$). Further, let $\overline V_B = \mu(B)V_B$, called, by analogy, {\it the normalized hitting time} 
(although the expected value of this variable need not be equal to 1). By ergodicity, $V_B$ and $\overline V_B$ are well defined. 
By an elementary consideration of the skyscraper above $B$, one easily verifies, that the distribution function 
$F_B$ of $\overline V_B$ satisfies, for every $t\ge 0$, the inequalities: 
$$
G_B(t) -\mu(B)\le F_B(t) \le G_B(t)
$$
(see [H-L-V] for more details). Because we deal with long blocks (so that, by the Shannon-McMillan-Breiman Theorem, $\mu(B)$ is, 
with high probability, very small), we will often replace $F_B$ by $G_B$. 

The key notions of this work are defined below:

\definition{Definition 1}
We say that the visits to $B$ {\it attract} (resp\. {\it repel}) each other with intensity $\epsilon$ from a 
distance $t>0$, if 
$$
F_B(t)\le 1-e^{-t}-\epsilon\ \ \  (\text{resp\. if }F_B(t)\ge 1-e^{-t}+\epsilon).
$$
We abbreviate that $B$ {\it attracts} ({\it repels}) with intensity $\epsilon$ if its visits attract (repel)
each other with intensity $\epsilon$ from some distance $t$.
\enddefinition

\definition{Definition 2}
We say that a process has {\it unbiased behavior} if there exist collections $\Cal B_n\subset \Cal P^n$ 
satisfying $\mu(\bigcup\Cal B_n)\to 1$, such that $F_{B_n}(t)\to 1-e^{-t}$ pointwise as $n\to\infty$, for any 
\sq\ of blocks $B_n\in\Cal B_n$.
\enddefinition

\definition{Definition 3}
We say that a process reveals {\it strong attracting}, if there is a subset $\na'\subset\na$ of upper density 1, 
and collections $\Cal B_{n'}\in\Cal P^{n'}$ for $n'\in\na'$, satisfying $\mu(\bigcup\Cal B_{n'})\to 1$, 
such that $F_{B_{n'}}(t)\to 0$ pointwise as $n'\to\infty$, for any \sq\ of blocks $B_{n'}\in\Cal B_{n'}$. 
\enddefinition

Let us explain why we use the terms ``attracting'' and ``repelling''. We will compare \proc\ 
with an independent Bernoulli process which is unbiased, i.e., for any long block $B$, $\tilde F_B(t) \approx 1-e^{-t}$ 
(and also $F_B(t) \approx 1-e^{-t}$) with high uniform accuracy (much better than $\epsilon$). Fix some $t>0$. Consider 
the random variable $I$ counting the number of occurrences of $B$ in the time period \hbox{$[0, \frac t{\mu(B)}]$}. 
The expected value of $I$ equals $\mu(B)\lfloor\frac t{\mu(B)}\rfloor\approx t$ (up to the ignorable error $\mu(B)$). 
On the other hand, $\mu\{I>0\} = \mu\{V_B\le \frac t{\mu(B)}\} = F_B(t)$. The ratio $\frac t{F_B(t)}$ represents
the conditional expected value of $I$ on the set $\{I>0\}$, i.e., the expected number of occurrences of $B$ in all 
intervals with at least one occurrence. Attracting from the distance $t$ means that $F_B(t)$ is smaller (by $\epsilon$) 
in \proc\ than in an independent Bernoulli process, i.e., that the above conditional expected value is larger in \proc\ 
than in the independent process. This fact can be further expressed as follows: If we observe the process \proc\ 
for time $\frac t{\mu(B)}$ (which is our ``memory length'' or ``lifetime of the observer'') and we happen to see the event 
$B$ during this time at least once, then the expected number of times we will observe the event $B$ is larger than the 
analogous value for a cylinder of the same measure in the independent Bernoulli process. The first occurrence of $B$ 
``attracts'' its further repetitions. The interpretation of repelling is symmetric. 

Obviously, occurrences of an event may simultaneously repel from one distance and attract from another. Notice, that the 
maximal intensity of repelling is $e^{-1}$ achieved at $t=1$ when $B$ appears periodically (this implies repelling from 
all distances). The intensity of attracting can be arbitrarily close to 1, which happens when $F_B(t)$ (hence also $G_B(t)$) 
remains near zero for some large $t$ (in particular this implies attracting from nearly all distances, except very small 
and very large ones, where marginal repelling can occur).
It is easy to see that such case happens exactly when the distribution of the normalized return time is nearly concentrated 
at zero, i.e., when most points in the set $B$ return after a time considerably smaller than $\frac 1{\mu(B)}$. Because the 
expected value of the return time equals $\frac 1{\mu(B)}$, there must be a small portion of $B$ with extremely large values 
of the return time. In such case the event $B$ appears in long series of high frequency, compensated by huge gaps of nearly 
complete absence. This is the essence of our notion of strong attracting. 

The first main result follows:

\proclaim{Theorem 1} If \proc\ is ergodic and has positive entropy, then for every \hbox{$\epsilon>0$} 
the measure of the union of all $n$-blocks $B\in\Cal P^n$ which repel with intensity $\epsilon$, converges 
to zero as $n$ grows to infinity.
\endproclaim


Obviously, Theorem~1 does not exclude the unbiased behavior. For example, a Bernoulli process with the independent generator
is unbiased. In fact, it follows from the results of [A-G], [H-S-V], that any process with a sufficient rate of mixing is 
unbiased (unbiased behavior is implied by ``exponential asymptotics''). Nevertheless, our second 
theorem will say in particular, that processes with the unbiased behavior are extremely exceptional among positive entropy processes.
\medskip
Let $(X,\mu)$ be a standard probability space, and let $m$ denote either a finite integer or the countable cardinal $\aleph_0$. 
The Rokhlin metric endows the collection of all measurable $\mu$-distinguishable partitions $\Cal P$ of $X$ into at most $m$ 
elements with a topology of a Polish space. 

\proclaim{Theorem 2} Let $(X,\mu,T)$ be an ergodic measure-preserving transformation of a standard probability space, with positive 
entropy. Fix some $2\le m \le \aleph_0$. Then, in the Polish space of all measurable partitions $\Cal P$ of $X$ into at most 
$m$ elements, there is a dense $G_\delta$ subset such that every partition in this subset generates a process which reveals strong attracting.
\endproclaim

Because partitions generating positive entropy form a dense open set (see Fact 5 below), 
we obtain that in a positive entropy measure preserving system a typical partition has both positive entropy and strong attracting.

\comment
\heading Interpretation and its limits\endheading

Strong attracting, i.e., attracting with intensity close to 1 occurs, when $G_B$ is very ``flat'' (close to zero on a long 
initial interval). Then $\tilde F_B$ is immediately very close to 1 indicating that on most of $B$ the first return time is 
much smaller than the mean $\frac 1{\mu(B)}$. Of course, this must be compensated on a small part of $B$ by extremely large values
of the return time. This means that the visits to $B$ occur in enormous clusters of very high frequency, compensated 
by huge pauses with no (or very few) visits.

\medskip\noindent
Repelling from the distance $t$ means exactly the opposite to attracting: The first occurrence lowers the expected number of 
repetitions within the observation period, i.e., repels them. If we have a mixed behavior, our impression about 
whether the event attracts or repels its repetitions depends on the length of our memory. Attracting not assisted 
by repelling (or assisted by repelling of an ignorably small intensity) means that no matter what memory length we apply, either we see a nearly unbiased behavior or the first occurrence visibly attracts further repetitions. 
Our Theorem~1 states that if we observe longer and longer blocks $B$, repelling from any distance must decay in both 
measure and intensity (while attracting can persist), so that for majority of long blocks we will see the behavior as 
described above.
We also note, that by pushing the graph of $F_B$ downward (compared to $1-e^{-t}$), attracting contributes 
to increasing the expected value of the associated random variable, i.e., of the hitting time. In case of 
attracting assisted by only very small intensity repelling, the average waiting time for the first occurrence 
of the event $B$ is increased in comparison to unbiased (may even not exist). Thus, instinctively judging the 
probability of the event by (the inverse of) the waiting time for the first occurrence we will typically 
underestimate it. All the more we are surprised, when the following occurrences happen after a considerably 
shorter time. This additionally strengthens the phenomenon's appearance.

Another consequence of attracting not assisted by repelling (or assisted by repelling of a very small intensity) 
is an increased variance of the return time statistic (the variance may even cease to exist). Thus, again, 
the gaps between the occurrences of $B$ are driven away from the expected value, toward the extremities 0
and $\infty$, and hence, into the pattern of clusters separated by longer pauses. We skip the elementary 
estimations of the variance. 

It must be reminded: Theorem 1 does not exclude that occurrences of long blocks will actually 
deviate from unbiased. There are conditions, weaker than full independence, under which the distributions of 
the normalized return times of long blocks converge almost surely to the exponential law. It is so, for instance,
in Markov processes (with finite memory). In fact, such convergence is implied by a sufficient rate of mixing 
([A-G], [H-S-V]). Yet, such processes seem to be somewhat exceptional and we expect that attracting rules in
majority of processes (see the Question 5 at the end of the paper). As we have already mentioned, at least that 
much is true, that in any dynamical system with positive entropy there exist partitions with strong attracting properties.
\endcomment

\heading More notation and preliminary facts\endheading
We now establish further notation and preliminaries needed in the proofs. If $\Bbb A\subset\z$ then we
will write $\Cal P^{\Bbb A}$ to denote the partition or sigma-field $\bigvee_{i\in \Bbb A}\sigma^{-i}(\Cal P)$.
We will abbreviate $\Cal P^n = \Cal P^{[0,n)}$, $\Cal P^{-n} = \Cal P^{[-n,-1]}$, $\Cal P^- = \Cal P^{(-\infty,-1]}$ 
(a ``finite future'', a ``finite past'', and the ``full past'' of the process).

We assume familiarity of the reader with the basics of entropy for finite partitions and sigma-fields
in a standard probability space. Our notation is compatible with [P] and we refer the reader to 
this book, as well as to [Sh] and [Wa], for background and proofs. In particular, we will be using the following:
\item * The entropy of a partition equals $H(\Cal P) = -\sum_{A\in\Cal P}\mu(A)\log_2(\mu(A))$.
\item * For two finite partitions $\Cal P$ and $\Cal B$, the conditional entropy $H(\Cal P|\Cal B)$ 
is equal to $\sum_{B\in\Cal B}\mu(B)H_B(\Cal P)$, where $H_B$ is the entropy evaluated for the conditional 
measure $\mu_B$ on $B$.
\item * The same formula holds for conditional entropy given a sub-sigma-field $\Cal C$, i.e., 
$$
\sum_{B\in\Cal B}\mu(B)H_B(\Cal P|\Cal C)=H(\Cal P|\Cal B\vee\Cal C). 
$$ 
\smallskip
\item * The entropy of the process is given by any one of the formulas below
$$
h = H(\Cal P|\Cal P^-) = \tfrac 1r H(\Cal P^r|\Cal P^-) = \lim_{r\to\infty}\tfrac 1r H(\Cal P^r).
$$

We will exploit the notion of $\epsilon$-independence for partitions and sigma-fields. The definition below
is an adaptation from [Sh], where it concerns finite partitions only. See also [Sm] for treatment of 
countable partitions. Because ``$\epsilon$'' is reserved for the intensity of repelling, we will speak 
about $\beta$-independence.

\definition{Definition 4} Fix $\beta>0$. 
A partition $\Cal P$ is said to be {\it $\beta$@-independent} of a sigma-field $\Cal B$ if for any
$\Cal B$-measurable countable partition $\Cal B'$ holds 
$$
\sum_{A\in\Cal P, B\in\Cal B'}|\mu(A\cap B)-\mu(A)\mu(B)|\le\beta.
$$ 
A process \proc\ is called a
{\it $\beta$-independent process} if $\Cal P$ is $\beta$-independent
of the past $\Cal P^-$.
\enddefinition

A partition $\Cal P$ is independent of another partition or a sigma-field $\Cal B$ if and only if 
$H(\Cal P|\Cal B) = H(\Cal P)$. The following approximate version of this fact holds (see [Sh, Lemma 7.3]
for finite partitions, from which the case of a sigma-field is easily derived). 

\proclaim{Fact 1}
A partition $\Cal P$ is $\beta$-independent of another partition or a sigma-field $\Cal B $ if 
$H(\Cal P|\Cal B)\ge H(\Cal P)-\xi$, for $\xi$ sufficiently small. \qed
\endproclaim
In course of the proof, a certain lengthy condition will be in frequent use. Let us introduce an abbreviation:

\definition{Definition 5}
Given a partition $\Cal P$ of a space with a probability measure $\mu$ and $\delta>0$, we will say that a 
property $\Phi(A)$ {\it holds for $A\in\Cal P$ with $\mu$-tolerance $\delta$} if 
$$
\mu\left(\bigcup\{A\in\Cal P: \Phi(A)\}\right)\ge 1-\delta. 
$$
\enddefinition

We shall also need an elementary estimate, whose proof is an easy exercise.

\proclaim{Fact 2}
For each $A\in\Cal P$, $H(\Cal P) \le (1-\mu(A))\log_2(\#\Cal P) + 1$. \qed
\endproclaim

In addition to the random variables of the
absolute and normalized return times $R_B$ and $\overline R_B$,
we will also use the analogous notions of the $k^{\text{th}}$ absolute return time
$$
R^{(k)}_B = \min\{i: \#\{0<j\le i: \sigma^j(y)\in B\} = k\},
$$ 
and of the normalized $k^{\text{th}}$ return time $\overline R^{(k)}_B = \mu(B)R^{(k)}_B$ (both defined on $B$), 
with $\tilde F^{(k)}_B$ always denoting the distribution function of the latter. Clearly, the expected value of 
$\overline R^{(k)}_B$ equals $k$. 

\heading The idea of the proof and the basic lemma\endheading
Before we pass to the formal proof of Theorem 1, we would like to have the reader oriented in the mainframe 
of the idea behind it. We intend to estimate (from above, by $1-e^{-t}+\epsilon$) the function $G_{B\!A}$
(replacing $F_{BA}$), for long blocks of the form $B\!A\in\Cal P^{[-n,r)}$.
The ``positive'' part $A$ has a fixed length $r$, while we allow the ``negative'' part $B$ to be arbitrarily long.
There are two key ingredients leading to the estimation. The first one, contained in Lemma 3, is the observation 
that for a fixed typical $B\in \Cal P^{-n}$, the part of the process induced on $B$ (with the conditional 
measure $\mu_B$) generated by the partition $\Cal P^r$, is not only a $\beta$-independent process, but it is 
also $\beta$-independent of many returns times $R_B^{(k)}$ of the cylinder $B$ (see the Figure 2).    

$$
\align
^{\text{coordinate }0}&_\downarrow \\ \vspace{-5pt}
...\boxed{\ \ _{_{B\phantom{_1}}} \ \ }\boxed{_{_{A_{\text{-}1}}}}...............\boxed{\ \ _{_{B\phantom{_1}}} \ \ }
&\boxed{_{_{A_0}}}..\boxed{\ \ _{_{B\phantom{_1}}} \ \ }\boxed{_{_{A_1}}}..........\boxed{\ \ _{_{B\phantom{_1}}} \ \ }\boxed{_{_{A_2}}}....\boxed{\ \
_{_{B\phantom{_1}}} \ \ }\boxed{_{_{A_3}}}....
\endalign
$$
{\eightpoint\it Figure 2: The process $\dots A_{-1}A_0A_1A_2\dots$ of $r$-blocks following the copies of 
$B$ is a $\beta$@-independent process with additional $\beta$-independence properties of 
the positioning of the copies of $B$.}

\medskip\noindent
This allows us to decompose (with high accuracy) the distribution function $\tilde F_{B\!A}$ of the normalized return 
time of $B\!A$ as follows:
$$
\gather
\tilde F_{B\!A}(t) = \mu_{B\!A}\{\overline R_{B\!A} \le t\} = \mu_{B\!A}\{R_{B\!A} \le \tfrac t{\mu(B\!A)}\} = \\
\sum_{k\ge 1} \mu_{B\!A}\{R^{(B)}_A = k, R^{(k)}_B \le \tfrac t{p\mu(B)}\}\approx 
\sum_{k\ge 1}\mu_{B\!A}\{R^{(B)}_A=k\}\cdot\mu_B\{\overline R^{(k)}_B \le \tfrac tp\} 
\approx \\ \sum_{k\ge 1}p(1-p)^{k-1}\cdot \tilde F^{(k)}_B(\tfrac tp),
\endgather
$$
where $R^{(B)}_A$ denotes the first (absolute) return time of $A$ in the process induced on $B$, and $p = \mu_B(A)$.

The second key observation is, assuming for simplicity full independence, that when trying to model some 
repelling for the blocks $B\!A$, we ascertain that it is largest, when the occurrences of $B$ are purely 
periodic. Any deviation from periodicity of the $B$'s may only lead to increasing the intensity of attracting 
between the copies of $B\!A$, never that of repelling. We will explain this phenomenon more formally in 
a moment. Now, if $B$ does appear periodically, then the normalized return time of $B\!A$ is governed by 
the same geometric distribution as the normalized return time of $A$ in the independent process induced
on $B$. If $p$ is small, this geometric distribution function becomes nearly the unbiased exponential 
law $1-e^{-t}$. The smallness of $p$ is {\it a priori} regulated by the choice of the parameter $r$ (Lemma~1). 

The phenomena that, assuming full independence, the repelling of $B\!A$ is maximized by periodic 
occurrences of $B$, and that even then there is nearly no repelling, is captured by the following 
elementary lemma, which will be also useful later, near the end of the rigorous proof.

\proclaim{Lemma 0}
Fix some $p\in(0,1)$. Let $\tilde F^{(k)}$ ($k\ge 1$) be a \sq\ of distribution functions on 
$[0,\infty)$ such that the expected value of the distribution associated to $\tilde F^{(k)}$ equals $k$. Define   
$$
\tilde F(t) = \sum_{k\ge 1}p(1-p)^{k-1} \tilde F^{(k)}(\tfrac tp),
\text{ \ \ and \ \ }
G(t) = \int_0^t 1-\tilde F(s)ds.
$$
Then $G(t)\le \frac1{\log e_p}(1-e_p^{-t})$, where $e_p = (1-p)^{-\frac 1p}$.
\endproclaim

\demo{Proof} We have
$$
G(t) = \sum_{k\ge 1}p(1-p)^{k-1}\int_0^t 1- \tilde F^{(k)}(\tfrac sp) ds.
$$
We know that $\tilde F^{(k)}(t)\in [0,1]$ and that $\int_0^\infty 1-\tilde F^{(k)}(s) ds = k$ (the expected value). 
With such constraints, it is the indicator function $1_{[k,\infty)}$ that maximizes the integrals from 
$0$ to $t$ simultaneously for every $t$ (because the ``mass'' $k$ above the graph is, for such choice 
of the function $\tilde F^{(k)}$, swept maximally to the left). The rest follows by direct calculations:
$$
\gather
G(t) \le \sum_{k\ge 1}p(1-p)^{k-1}\int_0^t 1_{[0,k)}(\tfrac sp) ds = 
\int_0^t \sum_{k=\lceil\frac sp\rceil}^\infty p(1-p)^{k-1} ds =\\ \int_0^t (1-p)^{\lceil\frac sp\rceil} ds \le \frac{(1-p)^{\frac tp}-1}{\log (1-p)^{\frac 1p}}.\qed
\endgather
$$ 
\enddemo

Recall that the maximizing distribution functions $\tilde F_B^{(k)} = 1_{[k,\infty)}$ occur, for the normalized 
return time of a set $B$, precisely when $B$ is visited periodically. This explains our former statement on
this subject.

\medskip
Let us comment a bit more on the first key ingredient, the $\beta$-independence. Establishing it is the most 
complicated part of the argument. The idea is to prove conditional (given a ``finite past'' $\Cal P^{-n}$) 
$\beta$-independence of the ``present'' $\Cal P^r$ from jointly the full past and a large part of the future,
responsible for the return times of majority of the blocks $B\in\Cal P^{-n}$. But the future part must not be 
too large. Let us mention the existence of ``bilaterally deterministic'' processes with positive entropy 
(first discovered by Gurevi\v c [G], see also [O-W1]), in which the sigma-fields generated by the coordinates 
$(-\infty,-m]\cup[m,\infty)$ do not decrease with $m$ to the Pinsker factor; they are all equal to the entire 
sigma-field. (Coincidently, our Example~1 has precisely this property; see the Remark 2.) Thus, in order 
to maintain any trace of independence of the ``present'' from our sigma-field
already containing the entire past, its part in the future must be selected with an extreme care. Let us also 
remark that an attempt to save on the future sigma-fields by adjusting them 
individually to each block $B_0\in\Cal P^{-n}$ falls short, mainly because of the ``off diagonal effect''; 
suppose $\Cal P^r$ is conditionally (given $\Cal P^{-n}$) nearly independent of a sigma-field which determines 
the return times of only one selected block $B_0\in \Cal P^{-n}$. The independence still holds conditionally 
given any cylinder $B\in\Cal P^{-n}$ from a collection of a large measure, but unfortunately, this collection 
can always miss the selected cylinder $B_0$. In Lemmas 2 and 3, we succeed in finding a sigma-field (containing 
the full past and a part of the future), of which $\Cal P^r$ is conditionally $\beta$-independent, and which  
``nearly determines'', for majority of blocks $B\in\Cal P^{-n}$, some finite number of their sequential return 
times (probably not all of them). This finite number is sufficient to allow the described earlier 
decomposition of the distribution function $\tilde F_{B\!A}$. 

\heading The proof of Theorem 1\endheading
Throughout the sequel we assume ergodicity and that the entropy $h$ of \proc\ is 
positive. We begin our computations with an auxiliary lemma allowing us to assume 
(by replacing $\Cal P$ by some $\Cal P^r$) that the elements of the ``present'' partition are small, 
relatively in most of $B\in\Cal P^n$ and for every $n$. Note that the Shannon-McMillan-Breiman Theorem
is insufficient: for the conditional measure the error term in that theorem depends increasingly on $n$, 
which we do not fix.

\proclaim{Lemma 1} For each $\delta$ there exists an $r\in\na$ such that 
for every $n\in\na$ the following holds for $B\in\Cal P^{-n}$ with $\mu$-tolerance $\delta$:
$$
\text{for every }A\in\Cal P^r,\ \ \mu_B(A)\le \delta.
$$
\endproclaim

\demo{Proof} 
Let $\alpha$ be so small that 
$$
\sqrt\alpha\le\delta \text{\ \ and\ \ } \frac{h-3\sqrt\alpha}{h+\alpha}\ge 1-\frac\delta2,
$$
and set $\gamma = \frac\alpha{\log_2(\#\Cal P)}$. Let $r$ be so big that 
$$
\frac 1r \le \alpha, \ \ \frac1{r(h+\alpha)}\le\frac\delta2,
$$
and that there exists a collection $\overline{\Cal P^r}$ of no more than $2^{r(h+\alpha)}-1$ elements of 
$\Cal P^r$ whose joint measure $\mu$ exceeds $1-\gamma$ (by the Shannon-McMillan-Breiman Theorem).

Let $\widetilde{\Cal P^r}$ denote the partition into the elements of $\overline{\Cal P^r}$ and the complement 
of their union, and let $\Cal R$ be the partition into the remaining elements of $\Cal P^r$ and the complement 
of their union, so that $\Cal P^r = \widetilde{\Cal P^r}\vee \Cal R$. For any $n$ we have 
$$
\gather
rh = H(\Cal P^r|\Cal P^-)\le H(\Cal P^r|\Cal P^{-n}) = 
H(\widetilde{\Cal P^r}\vee \Cal R|\Cal P^{-n})=\\
H(\widetilde{\Cal P^r}|\Cal R \vee \Cal P^{-n})+ H(\Cal R|\Cal P^{-n})\le
H(\widetilde{\Cal P^r}|\Cal P^{-n}) + H(\Cal R) \le  \\
\sum_{B\in\Cal P^{-n}}\mu(B)H_B(\widetilde{\Cal P^r}) + \gamma r\log_2(\#\Cal P) + 1
\endgather
$$
(we have used Fact 2 for the last passage). After dividing by $r$, we obtain
$$
\sum_{B\in\Cal P^{-n}}\mu(B)\tfrac 1r H_B(\widetilde{\Cal P^r}) \ge h - \gamma\log_2(\#\Cal P) - \tfrac1r 
\ge h - 2\alpha.
$$
Because each term $\tfrac 1r H_B(\widetilde{\Cal P^r})$ is not larger than $\frac 1r\log_2(\#\widetilde{\Cal P^r})$ 
which was set to be at most $h+\alpha$, we deduce that 
$$ 
\tfrac1r H_B(\widetilde{\Cal P^r}) \ge h - 3\sqrt\alpha
$$
holds for $B\in\Cal P^{-n}$ with $\mu$-tolerance $\sqrt\alpha$, hence also with $\mu$-tolerance $\delta$. 
On the other hand, by Fact 2, for any $B$ and $A\in\widetilde{\Cal P^r}$, holds: 
$$
H_B(\widetilde{\Cal P^r}) \le (1-\mu_B(A))\log_2(\#\widetilde{\Cal P^r})+1 \le (1-\mu_B(A))r(h+\alpha)+1.
$$
Combining the last two displayed inequalities we establish that, with $\mu$-tolerance $\delta$ for $B\in\Cal P^{-n}$ 
and then for every $A\in\widetilde{\Cal P^r}$, holds
$$
1-\mu_B(A)\ge \frac{h-3\sqrt\alpha}{h+\alpha} -\frac1{r(h+\alpha)}\ge 1-\delta.
$$
So, $\mu_B(A)\le \delta$. Because $\Cal P^r$ refines $\widetilde{\Cal P^r}$, the elements of $\Cal P^r$ 
are also not larger than~$\delta$.
\qed\enddemo

We continue the proof with a lemma which can be deduced from [Ru,~Lemma~3]. 
We provide a direct proof. For $\alpha>0$ and $M\in\na$ let 
$$
S(M,\alpha) = \bigcup_{m\in\z}[mM+\alpha M, (m+1)M-\alpha M)\cap\z.
$$

\proclaim{Lemma 2}
For fixed $\alpha$ and $r$ there exists $M_0$ such that for every $M\ge M_0$ holds, 
$$
H(\Cal P^r|\Cal P^-\vee \Cal P^{S(M,\alpha)})\ge rh-\alpha
$$ 
(see the Figure 3).
$$
***********\circ\circ..************..........************..........************..........
$$
\endproclaim\noindent
{\eightpoint\it Figure 3. The circles indicate the coordinates 0 through $r-1$, the conditioning sigma-filed 
is over the coordinates marked by stars, which includes the entire past and part of the future with gaps of size
$2\alpha M$ repeated periodically with period $M$ (the first gap is half the size).}

\demo{Proof} First assume that $r=1$. 
Denote also
$$
S'(M,\alpha) = \bigcup_{m\in\z}[mM+\alpha M, (m+1)M)\cap\z.
$$
Let $M$ be so large that $H(\Cal P^{(1-\alpha)M})< (1-\alpha)M(h+\gamma)$, where 
$\gamma = \frac{\alpha^2}{2(1-\alpha)}$. Then, for any $m\ge 1$,
$$
H(\Cal P^{S'(M,\alpha)\cap[0,mM)}|\Cal P^-)\le 
H(\Cal P^{S'(M,\alpha)\cap[0,mM)})<(1-\alpha)mM(h+\gamma).
$$ 
Because $H(\Cal P^{[0,mM)}|\Cal P^-) = mMh$, the complementary part of entropy must exceed 
$mMh - (1-\alpha)mM(h+\gamma)$ (which equals $\alpha mM (h-\tfrac\alpha2)$), i.e., we have
$$
H(\Cal P^{[0,mM)\setminus S'(M,\alpha)}|\Cal P^-\vee \Cal P^{S'(M,\alpha)\cap[0,mM)})
> \alpha mM (h-\tfrac\alpha2).
$$
Breaking the last entropy term as a sum over $j\in [0,mM)\setminus S'(M,\alpha)$ of the conditional
entropies of $\sigma^{-j}(\Cal P)$ given the sigma-field over all coordinates left of $j$ and all coordinates
from $S'(M,\alpha)\cap[0,mM)$ right of $j$, and because every such term is at most $h$, we deduce 
that more than half of these terms reach or exceed $h-\alpha$. So, a term not smaller than $h-\alpha$ 
occurs for a $j$ within one of the gaps in the left half of $[0,mM)$.
Shifting by $j$, we obtain $H(\Cal P|\Cal P^-\vee \sigma^i(\Cal P^{S'(M,\alpha)\cap [0,\frac{mM}2)}))\ge h-\alpha$,
where $i\in [0,\alpha M)$ denotes the relative position of $j$ in the gap. As we increase $m$, one value 
$i$ will repeat in this role along a subsequence $m'$. The operation $\vee$ is continuous for 
increasing \sq s of sigma-fields, hence $\Cal P^-\vee \sigma^i(\Cal P^{S'(M,\alpha)\cap [0,\frac{m'M}2)})$ 
converges over $m'$ to $\Cal P^-\vee \sigma^i(\Cal P^{S'(M,\alpha)})$. The entropy is continuous for such 
passage, hence $H(\Cal P|\Cal P^-\vee \sigma^i(P^{S'(M,\alpha)})\ge h-\alpha$. The assertion now follows 
because $S(M,\alpha)$ is contained in $S'(M,\alpha)$ shifted to the left by any $i\in[0,\alpha M)$.

Finally, if $r>1$, we can simply argue for $\Cal P^r$ replacing $\Cal P$. This will impose that
$M_0$ and $M$ are divisible by $r$, but it is not hard to see that for large $M$ the argument
works without divisibility at a cost of a slight adjustment of $\alpha$.
\qed\enddemo

For a long block $B\in\Cal P^{-n}$ let \procbr\ denote the process induced on $B$ generated by the restriction 
$\Cal P^r_B$ of $\Cal P^r$ to $B$ ($\sigma_B$ is the first return time map on $B$). 
The following lemma is the crucial item in our argument.

\proclaim{Lemma 3}
For every $\beta>0$, $r\in\na$ and $K\in\na$ there exists $n_0$ such that for every $n\ge n_0$, 
with $\mu$-tolerance $\beta$ for $B\in\Cal P^{-n}$, with respect to $\mu_B$, $\Cal P^r$ is $\beta$-independent 
of jointly the past $\Cal P^-$ and the first $K$ return times to $B$, $R^{(k)}_B$ ($k\in[1,K]$). 
In particular, \procbr\ is a $\beta$-independent process.
\endproclaim

\demo{Proof} 
We choose $\xi$ according to Fact 1, so that $\frac\beta2$-independence is implied.  Let $\alpha$ satisfy 
$$
0<\tfrac{2\alpha}{h-\alpha}<1, \ \ 18K\sqrt\alpha<1,\ \ \sqrt{2\alpha}<\xi,\ \ K\root 4\of\alpha<\tfrac\beta2.
$$ 
Let $n_0$ be so large that $H(\Cal P^r|\Cal P^{-n})<rh+\alpha$ for every $n\ge n_0$ and that for every 
$k\in [1,K]$ with $\mu$-tolerance $\alpha$ for $B\in\Cal P^{-n}$ holds
$$
\mu_B\{2^{n(h-\alpha)}\le R^{(k)}_B\le 2^{n(h+\alpha)}\}>1-\alpha
$$
(we are using Ornstein-Weiss Theorem [O-W2]; the multiplication by $k$, which should appear for the $k^{\text{th}}$ return time, 
is consumed by $\alpha$ in the exponent). 
Let $M_0\ge 2^{n_0(h-\alpha)}$ be so large that the assertion of Lemma~2 holds for $\alpha$, $r$ and $M_0$, 
and that for every $M\ge M_0$, 
$$
(M+1)^{1+\frac{2\alpha}{h-\alpha}}<\alpha M^2 \text{ \ and \ }\tfrac{\log_2(M+1)}{M(h-\alpha)}<\alpha. 
$$
We can now redefine (enlarge) $n_0$ and $M_0$ so that $M_0=\lfloor 2^{n_0(h-\alpha)}\rfloor$. Similarly,
for each $n\ge n_0$ we set $M_n=\lfloor 2^{n(h-\alpha)}\rfloor$. Observe, that the interval where the first 
$K$ returns of most $n$-blocks $B$ may occur (up to probability $\alpha$), is contained in $[M_n, \alpha M_n^2]$
(because $2^{n(h+\alpha)}\ \le (M_n+1)^{1+\frac{2\alpha}{h-\alpha}}<\alpha M_n^2$). 

At this point we fix some 
$n\ge n_0$. The idea is to carefully select an $M$ between $M_n$ and $2M_n$ (hence not smaller than $M_0$), 
such that the initial $K$ returns of nearly every $n$-block happen most likely inside (with all its $n$ 
symbols) the set $S(M,\alpha)$, so that they are ``controlled'' by the sigma-field $\Cal P^{S(M,\alpha)}$. 
Let $\alpha' = \alpha + \frac n{M_n}$, so that every $n$-block overlapping with $S(M,\alpha')$ is completely
covered by $S(M,\alpha)$. By the second assumption on $M\ge M_0$ and by the formula connecting $M_n$ and $n$, 
we have $\alpha'<2\alpha$. To define $M$ we will invoke the triple Fubini Theorem. Fix $k\in[1,K]$ and consider 
the probability space 
$$
\Cal P^{-n}\times [M_n,2M_n]\times \na
$$ 
equipped with the (discrete) measure $\Cal M$ whose marginal on $\Cal P^{-n}\times [M_n,2M_n]$ is the product 
of $\mu$ (more precisely, of its projection onto $\Cal P^{-n}$) with the uniform distribution on the integers 
in $[M_n,2M_n]$, while, for fixed $B$ and $M$, the measure on the corresponding $\na$-section is the distribution 
of the random variable $R_B^{(k)}$. In this space let $S$ be the set whose $\na$-section for a fixed 
$M$ (and any fixed $B$) is the set $S(M,\alpha')$. We claim that for every $l\in [M_n,\alpha M_n^2]\cap\na$ 
(and any fixed $B$) the $[M_n,2M_n]$-section of $S$ has measure exceeding $1-16\alpha$. This is quite obvious 
(even for every $l\in [M_n,\infty)$ and with $1-15\alpha$) if $[M_n,2M_n]$ is equipped with the normalized 
Lebesgue measure (see the Figure 4). 

\input epsf.tex
\epsfxsize=13truecm
\centerline{\epsfbox{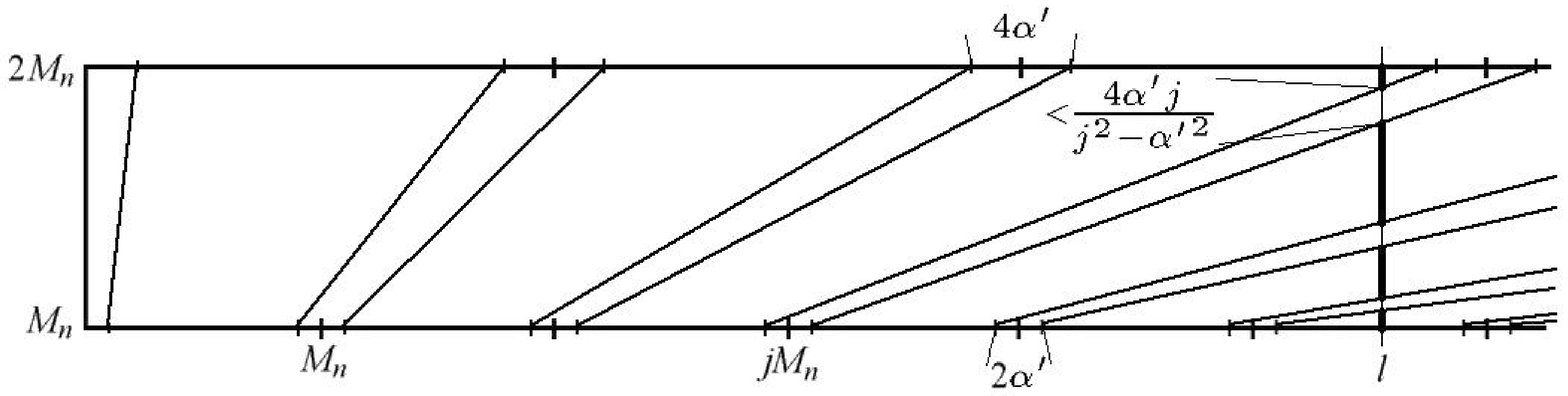}}
\noindent{\eightpoint\it Figure 4: The complement of $S$ splits into thin skew strips shown in the picture.
The normalized Lebesgue measure of any vertical section of the $j^{\text{\rm th}}$ strip (starting at $jM_n$ 
with $j\ge 1$) is at most $\frac{4\alpha'j}{j^2 - {\alpha'}^2}\le \frac{5\alpha'}j\le \frac{10\alpha}j$. 
Each vertical line at $l\ge M_n$ intersects strips with indices $j, j+1, j+2$ up to at most $2j$ 
(for some $j$), so the joint measure of the complement of the section of $S$ does not exceed $15\alpha$.}

\medskip
$$
\align \ \ 
&\ \ \ \ \ \ \ \ \ \ \ \ \ \ \ \ \ \ S\\
&\ \ \ \ \ \ \ \ \ \ \ \ \swarrow\ \ \downarrow\ \ \ \searrow\\
\vspace{-5pt}
&|\!\overline{..}.\overline{...}.\overline{....}.\overline{.....}..\overline{......}..
\overline{.......}..\overline{........}...\overline{.........}...\overline{..........}...\overline{..........}...\overline{.....}\!| \\
&\!M_n\hphantom{...................................................................................}2M_n
\endalign
$$
{\eightpoint\it Figure 5: The discretization replaces the Lebesgue measure by the uniform measure on 
$M_n$ integers, thus the measure of any interval can deviate from its Lebesgue measure by at most $\frac 1{M_n}$.
For $l\le \alpha M_n^2$ the corresponding section of $S$ (in this picture drawn horizontally) 
consists of at most $\alpha M_n$ intervals, so its measure can deviate by no more than $\alpha$.}
\medskip
In the discrete case, however, {\it a priori} it might happen that the integers
along some $[M_n,2M_n]$-section often ``miss'' the section of $S$ leading to a decreased measure value. (For 
example, it is easy to see that for $l = (2M_n)!$ the measure of the section of $S$ is zero.) But because we 
restrict to $l\le \alpha M_n^2$, the discretization does not affect the measure of the section of $S$ 
by more than $\alpha$, and the estimate with $1-16\alpha$ holds (see the Figure 5 above).

Taking into account all other inaccuracies (the smaller than $\alpha$ part of $S$ outside $[M_n,\alpha M_n^2]$ 
and the smaller than $\alpha$ part of $S$ projecting onto blocks $B$ which do not obey the Ornstein-Weiss return 
time estimate) it is safe to claim that 
$$
\Cal M(S)>1-18\alpha.
$$ 
This implies that for every $M$ from a set of measure at least $1-18\sqrt{\alpha}$ 
the measure of the $(\Cal P^{-n}\times \na)$-section of $S$ is larger than or equal to $1-\sqrt\alpha$. 
For every such $M$, with $\mu$-tolerance $\root4\of\alpha$ for $B\in\Cal P^{-n}$, the probability $\mu_B$ 
that the $k^{\text{\rm th}}$ repetition of $B$ falls in $S(M,\alpha')$ (hence with all its $n$ terms inside 
the set $S(M,\alpha)$) is at least $1-\root 4\of\alpha$.

Because $18K\sqrt\alpha<1$, there exists at least one $M$ for which the above holds for every $k\in[1,K]$.
This is our final choice of $M$ which from now on remains fixed. For this $M$, and for cylinders $B$ chosen 
with $\mu$-tolerance $K\root 4\of\alpha$, each of the considered $K$ returns of $B$ with probability 
$1-\root 4\of\alpha$ falls (with all its coordinates) inside $S(M,\alpha)$. Thus, for such a $B$, with probability 
$1-K\root 4\of\alpha$ the same holds simultaneously for all $K$ return times. In other words, 
there is a set $U_{B}$ of measure not exceeding $K\root 4\of\alpha$ outside of which $R^{(k)}_B = \tilde R^{(k)}_B$, 
where $\tilde R^{(k)}_B$ is defined as the time of the $k^{\text{th}}$ fully visible inside $S(M,\alpha)$ return of $B$.
Notice that $\tilde R^{(k)}_B$ is $\Cal P^{S(M,\alpha)}$-measurable.

Let us go back to our entropy estimates. We have, by Lemma 2,
$$
\gather
\sum_{B\in\Cal P^{-n}}\mu(B)H_B(\Cal P^r|\Cal P^-\vee \Cal P^{S(M,\alpha)}) = 
H(\Cal P^r|\Cal P^{-n}\vee\Cal P^-\vee \Cal P^{S(M,\alpha)})=\\
H(\Cal P^r|\Cal P^-\vee \Cal P^{S(M,\alpha)}) \ge
rh -\alpha \ge H(\Cal P^r|\Cal P^{-n})-2\alpha =\\
\sum_{B\in\Cal P^{-n}}\mu(B)H_B(\Cal P^r)-2\alpha.
\endgather
$$
Because $H_B(\Cal P^r|\Cal P^-\vee \Cal P^{S(M,\alpha)})\le H_B(\Cal P^r)$ for every $B$,
we deduce that with $\mu$-tolerance $\sqrt{2\alpha}$ for $B\in\Cal P^{-n}$ must hold 
$$
H_B(\Cal P^r|\Cal P^-\vee \Cal P^{S(M,\alpha)})\ge H_B(\Cal P^r)-\sqrt{2\alpha}\ge H_B(\Cal P^r)-\xi.
$$
Combining this with the preceding arguments, with $\mu$-tolerance $K\root 4\of\alpha + \sqrt{2\alpha}<\beta$
for $B\in\Cal P^{-n}$ both the above entropy inequality holds, and we have the estimates of the measures of sets $U_B$. 
By the choice of $\xi$, we obtain that with respect to $\mu_B$, $\Cal P^r$ is jointly $\frac\beta2$-independent 
of the past and the modified return times $\tilde R^{(k)}_B$ ($k\in [1,K]$). Because $\mu(U_B)\le K\root4\of\alpha<\frac\beta2$, 
this clearly implies $\beta$-independence if each $\tilde R^{(k)}_B$ is replaced by $R^{(k)}_B$.
\qed\enddemo

To complete the proof of Theorem 1 it now remains to put the items together.

\demo{Proof of Theorem 1}
Fix an $\epsilon>0$. On $[0,\infty)$, the functions 
$$
g_p(t) = \min\{1,\tfrac1{\log e_p}(1-e_p^{-t}) +pt\},
$$ 
where $e_p = (1-p)^{-\frac 1p}$, decrease uniformly to $1-e^{-t}$ as $p\to 0^+$. So, let $\delta$ be such 
that $g_{\delta}(t)\le 1-e^{-t}+\epsilon$ for every $t$. We also assume that 
$$
(1-2\delta)(1-\delta)\ge 1-\epsilon.
$$
Let $r$ be specified by Lemma 1, so that $\mu_B(A)\le \delta$ for every $n\ge 1$, every $A\in\Cal P^r$
and for $B\in\Cal P^{-n}$ with $\mu$-tolerance $\delta$. On the other hand, once $r$ is fixed, the 
partition $\Cal P^r$ has at most 
$(\#\Cal P)^r$ elements, so with $\mu_B$-tolerance $\delta$ for $A\in\Cal P^r$, $\mu_B(A)\ge \delta(\#\Cal P)^{-r}$. 
Let $\Cal A_B$ be the subfamily of $\Cal P^r$ (depending on $B$) where this inequality holds.
Let $K$ be so large that for any $p\ge \delta(\#\Cal P)^{-r}$,
$$
\sum_{k=K+1}^\infty p(1-p)^k <\tfrac\delta2,
$$
and choose $\beta<\delta$ so small that 
$$
(K^2+K+1)\beta<\tfrac\delta2.
$$ 
The application of Lemma 3 now provides an $n_0$ such that for any $n\ge n_0$, with $\mu$-tolerance $\beta$
for $B\in\Cal P^{-n}$, the process induced on $B$ generated by $\Cal P^r$ has the desired $\beta$-independence 
properties involving the initial $K$ return times of $B$. So, with tolerance $\delta+\beta<2\delta$ we have 
both, the above $\beta$-independence and the estimate $\mu_B(A)<\delta$ for every $A\in\Cal P^r$. Let $\Cal B_n$ 
be the subfamily of $\Cal P^{-n}$ where these two conditions hold. Fix some $n\ge n_0$.

Let us consider a cylinder set $B\cap A\in\Cal P^{[-n,r)}$ (or, equivalently, the block $B\!A$), 
where $B\in\Cal B_n$, $A\in \Cal A_B$. The length of $B\!A$ is $n+r$, which represents an arbitrary integer 
larger than $n_0+r$. Notice that the family of such sets $B\!A$ covers more than $(1-2\delta)(1-\delta)\ge 1-\epsilon$ 
of the space. 

We will examine the distribution of the normalized first return time for $B\!A$. 
In addition to our customary notations of return times, let $R^{(B)}_A$ be the first (absolute) 
return time of $A$ in \procbr, i.e., the variable defined on $B\!A$, counting the number of visits to 
$B$ until the first return to $B\!A$. Let $p = \mu_B(A)$ (recall, this is not smaller than $\delta(\#\Cal P)^{-r}$). 
We have
$$
\gather
\tilde F_{B\!A}(t) = \mu_{B\!A}\{\overline R_{B\!A} \le t\} = \mu_{B\!A}\{R_{B\!A} \le \tfrac t{\mu(B\!A)}\} = \\
\sum_{k\ge 1} \mu_{B\!A}\{R^{(B)}_A = k, R_B^{(k)} \le \tfrac t{p\mu(B)}\}.
\endgather
$$ 
The $k^{\text{th}}$ term of this sum equals
$$
\tfrac 1p \mu_B(\{A_k = A\}\cap\{A_{k-1}\neq A\}\cap\dots
\cap\{A_1\neq A\}\cap\{A_0 = A\}\cap\{R^{(k)}_B \le \tfrac t{p\mu(B)}\}),
$$
where $A_i$ is the $r$-block following the $i^{\text{th}}$ copy of $B$ (the counting starts from $0$ at
the copy of $B$ positioned at $[-n,-1]$).

By Lemma 3, for $k\le K$, in this intersection of sets each term is $\beta$-independent of the intersection right from it. So, proceeding from the left, we can replace the probabilities of the
intersections  by products of probabilities, allowing an error of $\beta$. Note that the last term equals 
$\mu_B\{\overline R^{(k)}_B \le \tfrac tp\} = \tilde F^{(k)}_B(\tfrac tp)$. Jointly, the inaccuracy will not exceed 
$(K+1)\beta$:
$$
\left|\mu_{B\!A}\{R^{(B)}_A = k, R_B^{(k)} \le \tfrac t{p\mu(B)}\} - p(1-p)^{k-1}\tilde F^{(k)}_B(\tfrac tp)\right|
\le (K+1)\beta.
$$
Similarly, we also have $\left|\mu_{B\!A}\{R^{(B)}_A=k\}-p(1-p)^{k-1}\right|\le K\beta$, hence the tail of the series 
$\mu_{B\!A}\{R^{(B)}_A=k\}$ above $K$ is smaller than $K^2\beta$ plus the tail of the geometric series $p(1-p)^{k-1}$,
which, by the fact that $p\ge \delta(\#\Cal P)^{-r}$, is smaller than $\tfrac\delta2$. Therefore
$$
\tilde F_{B\!A}(t) \approx \sum_{k\ge 1}p(1-p)^{k-1}\tilde F^{(k)}_B(\tfrac tp),
$$
up to $(K^2+K+1)\beta +\tfrac\delta2\le\delta$, uniformly for every $t$. By the application of Lemma~0, 
$G_{B\!A}$ satisfies
$$
G_{B\!A}(t) \le \min\{1, \tfrac1{\log e_p}(1-e_p^{-t}) + \delta t\}\le 
g_{\delta}(t) \le 1-e^t+\epsilon
$$
(because $p\le\delta$).
We have proved that for our choice of $\epsilon$ and an arbitrary length $m\ge n_0+r$, with $\mu$-tolerance 
$\epsilon$ for the cylinders $C\in \Cal P^m$, the intensity of repelling between visits to $C$ is at most 
$\epsilon$. This concludes the proof of Theorem 1.
\qed\enddemo

\heading Proof of Theorem 2 \endheading

This proof requires a number of technical ingredients, such as ``semi-periodic markers'' or short ``transciently forbidden words''. 
The two facts below are standard exercises in ergodic theory and we only outline their proofs.

\proclaim{Fact 3}
In a process \proc\ of positive entropy, where $\Cal P$ is finite or countable, for each $k\in\na$ and $\epsilon>0$ 
there exist an $l\in\na$ and $k$ words $w_1, w_2,\dots, w_k$ of length~$l$ such that \newline
1. each $w_i$ starts and ends with the same symbol $a\in\Cal P$, independent from $i$\newline
2. each $w_i$ has measure $\mu$ at most $\frac\epsilon{lk}$,\newline
3. for each $i$ the set 
$$
w_i\setminus \bigcup_{j\neq i}\bigcup_{m=-l}^l \sigma^m(w_j)
$$
has positive measure $\mu$.

\endproclaim
\demo{Proof} For 1\. use recurrence in the $k$-fold product system, and for 2\. use the Shannon-McMillan-Breiman Theorem. Condition 3\. follows easily from the high complexity in positive entropy. 
\qed\enddemo

\proclaim{Fact 4} In every measure-preserving system $(X,\mu,T)$ of positive entropy $h$, for each sufficiently large $r\in\na$ 
there exists a ``semiperiodic $r$-marker'', i.e., a measurable set $F$ such that the first return time $R_F$ assumes only two values: $r$ and $r+1$. 
\endproclaim
\demo{Proof} The system has a Bernoulli factor of entropy $h$. For large $r$ the binary process obtained by random concatenations 
of two blocks, $0^{r-1}1$ and $0^r1$, is Bernoulli with entropy smaller than $h$, hence it is a factor of $(X,\mu,T)$. The
lift of the cylinder over $1$ is the desired set $F$ in $X$.
\qed\enddemo

We are in a position to present the proof of Theorem 2.

\demo{Proof of Theorem 2} 
Fix $\epsilon>0$, $t>0$ and $N\in\na$. Consider the following property of a (finite or countable) partition $\Cal P$: for every 
$n\in[N,N^2]$, $F_B(t)<\epsilon$ with $\mu$-tolerance $\epsilon$ for $B\in\Cal P^n$. (Recall that $F_B$ denotes the distribution
function of the normalized hitting time for $B$). It is easy to see that it holds on an open set $\Cal E_{\epsilon, t, N}$ of 
partitions (both in the space of partitions into at most $m$ 
elements and in the space of at most countable partitions); for each $n$ we can take the same finite sets of ``good'' $n$-cylinders 
$B$ for the partitions in a neighborhood of $\Cal P$ as for $\Cal P$. Of course, the set 
$$
\Cal E_{\epsilon, t} = \bigcup_{N\ge 1}\Cal E_{\epsilon, t, N},
$$ 
of partitions such that the same property holds for some $N$, is also open. The main effort in the proof will be to show 
that this set is also dense. Once this is done, the proof is complete, because then the dense $G_\delta$ set of partitions 
which reveal strong attracting can be obtained by intersecting the sets $\Cal E_{\epsilon, t}$ over countably many pairs 
$(\epsilon, t)$ with $\epsilon\to 0$ and $t\to\infty$. Notice that for any infinite \sq\ of natural numbers $N$ the set $\bigcup[N,N^2]$ has upper density 1 in $\na$.

In order to prove the density of $\Cal E_{\epsilon, t}$, fix a (finite or countable) partition $\Cal P$. Set
$$
k = \lceil \tfrac2\epsilon\rceil+1, \ \ \delta=\tfrac {1-\frac\epsilon2}{4k}, \ \ M=2kt.
$$

Choose words
$w_1,w_2,\dots, w_k$ according to Fact 3. Let $N$ be so large, that with $\mu$-tolerance $\frac\epsilon2$ in every $N$-block,
every word $w_i$ occurs at least once so that it does not overlap with any other $w_j$ (see condition 3. in Fact 3). 
Obviously, the same holds if $N$ is replaced by any larger integer. For every $n\in [N,N^2]$ we can thus select a finite 
collection of ``good'' $n$-blocks which satisfy the above and cover $1-\frac\epsilon2$ of the space. Let $p$ be so large, that
$\frac{2N^2}p<\frac\epsilon2$, and that every good $n$-block (for any $n\in [N,N^2]$) occurs at least $M$ times in every, 
up to $\mu$-tolerance $\delta$, $\frac p2$-block. Let $r = kp$.  

Now we invoke the semiperiodic $r$-marker set $F$ of Fact 4. Every $\Cal P$-name can be divided at visits to $F$ into a 
concatenation of $r$-blocks and $(r+1)$-blocks. For simplicity, we will call all of them {\it component $r$-blocks}. 
Every component $r$-block $C$ will be further decomposed as a concatenation of $k$ 
$p$-blocks $C_1C_2\dots C_k$ ($C_k$ is either a $p$-block or a $(p+1)$-block, but again, for simplicity, we will cal all 
these blocks $p$-blocks). We fix a symbol $b\neq a$ in $\Cal P$ (recall that $a$ denotes the first and last symbol of each $w_i$). 
Now we modify the partition $\Cal P$ by changing the $\Cal P$-names of points, as follows: In every $\Cal P$-name we replace, for every $i$, every occurrence of $w_i$ within every $i^{\text{th}}$ $p$-block $C_i$ of every component $r$-block $C$ and within the
first $N^2$ positions of the following $p$-block $C_{i+1}$ (here $k+1 = 1$), by the word $w_0=b^l$. Notice that there is no 
collision when overlapping words are replaced.
$$
\align
...|&
\underset{C_1}\to{\underbrace{\overset{N^2}\to{\overbrace{.w_2....}}.w_3.w_1.w_1...w_2...w_3.w_1...w_2.}}|
\underset{C_2}\to{\underbrace{\overset{N^2}\to{\overbrace{.w_2.w_1}}...w_3...w_2..w_3....w_1.w_2.w_2...}}|
\overset{N^2}\to{\overbrace{.w_1.w_2}}...\\
...|&
\underset{C_1}\to{\underbrace{\overset{N^2}\to{\overbrace{.w_2....}}.w_3.w_0.w_0...w_2...w_3.w_0...w_2.}}|
\underset{C_2}\to{\underbrace{\overset{N^2}\to{\overbrace{.w_0.w_0}}...w_3...w_0..w_3....w_1.w_0.w_0...}}|
\overset{N^2}\to{\overbrace{.w_1.w_0}}...\\
\endalign
$$
{\eightpoint\it Figure 6: A  $\Cal P$-name before and after modification.}
\bigskip\noindent

Let $C_i'$ denote the right part of $C_i$ obtained by cutting off its left $N^2$ entries.

First observe, that the change affects only a subset of 
$\bigcup_{i=1}^k\bigcup_{m=0}^{l-1} \sigma^{-m}(w_i)$, whose measure is smaller than $\epsilon$. Thus the distance between
$\Cal P$ and the partition $\Cal P'$ after the modification is less than $\epsilon$.

Notice also, that the modification completely forbids the word $w_i$ within any $C_i$ and $N^2$ positions right from it, 
because all ``old'' occurrences are removed, and the insertions of the block $w_0$ do not create any overlapping 
``new'' instances of $w_i$. On the other hand, these modifications do not affect inside $C'_i$ the words $w_j$ with 
$j\neq i$ which have not overlapped with $w_i$ before the change.

For fixed $n\in[N,N^2]$ and $i\in [1,k]$ observe an $n$-block $B'$ (over the partition $\Cal P'$) obtained from a 
``good'' $n$-block $B$ over $\Cal P$ appearing inside some $C'_i$. Such blocks (with all possible values of $i$) still 
cover more than $1-\frac\epsilon2-\frac{2N^2}p\ge 1-\epsilon$ of the space.
Because, for each $j\neq i$, $B'$ contains at least one unaffected copy of $w_j$ (not overlapping with $w_i$ before the
change), $B'$ cannot occur with its leftmost position located in any $C_j$ except for $j=i$. On the other hand, inside  
$C'_i$ it occurs as many times as $B$ did before the change. Because the blocks $C'_i$ jointly contain a fraction
$\frac{\frac p2-N^2}{kp} \ge \frac {1-\frac\epsilon2}{2k}$ of all $\frac p2$-blocks, only a fraction of at most 
$\frac{2k\delta}{1-\frac\epsilon2} = \frac12$ of all blocks $C'_i$ may contain less than $M$ copies of $B'$. Thus the 
measure $\mu(B')$ of the cylinder $B'$ (with respect to the partition $\Cal P'$) is at least $\frac M{2kp} = \frac tp$.
The waiting time for $B'$ is not larger than $p$ only within $C_i$ and the preceding $p$-block, so 
$\mu\{V_{B'}\le p\}\le\frac 2k < \epsilon$. After normalizing, we obtain $F_{B'}(t)<\epsilon$. We have proved that 
$\Cal P' \in \Cal E_{\epsilon,t,N}$. This completes the proof of the claim that $\Cal E_{\epsilon,t}$ is dense among 
the partitions, and ends the whole proof.
\qed\enddemo

For a more complete image of a process generated by a typical partition, let us formulate one more fact.

\proclaim{Fact 5}
Let $(X,\mu,T)$ be an ergodic measure-preserving transformation of with positive entropy. Fix some $2\le m \le \aleph_0$. 
Then, in the Polsh space of all measurable partitions $\Cal P$ of $X$ into at most $m$ elements, the set of
partitions generating positive entropy is open and dense.
\endproclaim
\demo{Proof}
It is known that entropy is continuous in the Rokhlin metric, so positive entropy is an open property (see e.g. [P]). 
To obtain density it suffices to perturb a zero-entropy partition by a small set not measurable with respect to
the Pinsker algebra.\qed\enddemo

\heading Consequences for limit laws\endheading
The studies of limit laws for return/hitting time statistics are based on the following approach:
For $x\in\Cal P^\z$ define $F_{x,n} = F_B$ (and $\tilde F_{x,n}=\tilde F_B$), where $B$ is the block $x[0,n)$ 
(or the cylinder in $\Cal P^n$ containing $x$). Because for nondecreasing functions $F:[0,\infty)\to[0,1]$, the 
weak convergence coincides with the convergence at continuity points, and it makes the space of such functions 
metric and compact, for every $x$ there exists a well defined collection of limit distributions for $F_{x,n}$ 
(and for $\tilde F_{x,n}$) as $n\to \infty$. They are called {\it limit laws for the hitting (return) 
times at $x$}. Due to the integral relation ($F_B\approx G_B$) a \sq\ of return time distributions 
converges weakly if and only if the corresponding hitting time distributions converge pointwise (see [H-L-V]), 
so the limit laws for the return times completely determine those for hitting times and {\it vice versa}. 
A limit law is {\it essential} \,if it appears along some sub\sq\ $(n_k)$ for $x$'s in a set of positive measure. 
In particular, the strongest situation occurs when there exists an almost sure limit law along the full \sq\ $(n)$.
In such case the process is said to have {\it exponential asymptotics}.
Most of the results concerning the limit laws, obtained so far, can be classified in three major groups:
\item\item{a)} characterizations of possible essential limit laws for specific zero entropy processes (e.g. [D-M], [C-K];
these limit laws are usually atomic for return times or piecewise linear for hitting times), 
\item\item{b)} finding classes of processes with exponential asymptotics (e.g. [A-G], [H-S-V]), and 
\item\item{c)} results concerning 
non-essential limit laws, limit laws along sets other than cylinders (see [L]; every probabilistic distribution 
with expected value not exceeding 1 can occur in any process as the limit law for such general return times),  
or other very specific topics. 

\noindent As a consequence of our Theorem 1, we obtain, for the first time, a serious 
bound on the possible essential limit laws for the hitting time statistics along cylinders in the 
general class of ergodic positive entropy processes. The statement (1) below is even slightly stronger, 
because we require, for a subsequence, convergence on a positive measure set, but not necessarily to a common limit.

\proclaim{Theorem 3} Assume ergodicity and positive entropy of the process \proc.
\roster
\item If a sub\sq\ $(n_k)$ is such that
$\tilde F_{x,n_k}$ converge pointwise to some limit laws $\tilde F_x$ on a positive measure set 
$A$ of points $x$, then almost surely on $A$, $\tilde F_x(t)\le 1-e^{-t}$ at each $t\ge0$. 
\item If $(n_k)$ grows sufficiently fast, then there is a full measure set, such that for every $x$
in this set holds: $\limsup_k\, \tilde F_{x,n_k}(t)\le 1-e^{-t}$ at each $t\ge0$.
\endroster
\endproclaim

\demo{Proof}
The implication from Theorem 1 to Theorem 3 is obvious and we leave it to the reader. 
For (2) we hint that $(n_k)$ must grow fast enough to ensure summability of the measures of the sets 
where the intensity of repelling persists, then the Borel-Cantelli Lemma applies.
\qed\enddemo

Our Theorem 2 (again combined with the Borel-Cantelli Lemma) shows that a typical positive entropy process 
(including Bernoulli processes) admits the zero function as an essential limit law for the distributions 
of the waiting time. In particular, not all Bernoulli processes have exponential asymptotics. 

\heading An example \endheading

It is important not to be misled by an oversimplified approach to Theorem~1. The ``decay of repelling'' in positive 
entropy processes appears to agree with the intuitive understanding of entropy as chaos: repelling is a ``self-organizing'' 
property; it leads to a more uniform, hence less chaotic, distribution of an event along a typical orbit. Thus 
one might expect that repelling with intensity $\epsilon$ revealed by a fraction $\xi$ of all $n$-blocks 
contributes to lowering an upper estimate of the entropy by some percentage proportional to $\xi$ and depending 
increasingly on $\epsilon$. If this happens for infinitely many lengths $n$ with the same parameters $\xi$ and 
$\epsilon$, the entropy should be driven to zero by a geometric progression. Surprisingly, it is not quite so, 
and the phenomenon has more subtle grounds. We will present an example which exhibits the incorrectness of 
such intuition. Note also that in the proof of Theorem~1 the entropy is ``killed completely in one step'', that means, 
positive entropy and persistent repelling lead to a contradiction by examining the blocks of one sufficiently large length 
$n$; we do not use any iterated procedure requiring repelling for infinitely many lengths.

The construction below will show that for each $\delta>0$ and $n\in\na$ there exists $N\in\na$ and an ergodic 
process on $N$ symbols with entropy $\log_2N -\delta$, such that the $n$-blocks from a collection of joint 
measure equal to $\frac 1n$ repel with nearly the maximal possible intensity $e^{-1}$. Because $\delta$ 
can be extremely small compared to $\frac 1n$, this construction illustrates, that there is no ``reduction of entropy'' by an amount proportional to the fraction of blocks which reveal strong repelling. 

\example{Example 1} Let $\Cal P$ be an alphabet 
of a large cardinality $N$. Divide $\Cal P$ into two disjoint subsets, one, denoted $\Cal P_0$, of cardinality 
$N_0 = N2^{-\delta}$ and the relatively small (but still very large) rest which we denote by $\{1,2,\dots,r\}$ 
(we will refer to these symbols as ``markers''). For $i = 1,2,\dots,r$, let $\Cal B_i$ be the collection of all 
$n$-blocks whose first $n-1$ symbols belong to $\Cal P_0$ and the terminal symbol is the marker $i$. The cardinality 
of $\Cal B_i$ is $N_0^{n-1}$. Let $\Cal C_i$ be the collection of all blocks of length $nN_0^{n-1}$ obtained 
as concatenations of blocks from $\Cal B_i$ using each of them exactly once. The cardinality of $\Cal C_i$ 
is $(N_0^{n-1})!$. Let $X$ be the subshift whose points are infinite concatenations of blocks from 
$\bigcup_{i=1}^r\Cal C_i$, in which every block belonging to $\Cal C_i$ is followed by a block from 
$\Cal C_{i+1}$ ($1\le i<r$) and every block belonging to $\Cal C_r$ is followed by a block from $\Cal C_1$.
Let $\mu$ be the shift-invariant measure of maximal entropy on $X$. It is immediate to see that the entropy
of $\mu$ is $\frac 1{nN_0^{n-1}}\log_2((N_0^{n-1})!)$, which, for large $N$, nearly equals $\log_2{N_0} = \log_2N -\delta$. Finally observe that the measure of each $B\in\Cal B_i$ equals $\frac 1{nrN_0^{n-1}}$, the joint measure 
of $\bigcup_{i=1}^r\Cal B_i$ is exactly $\frac1n$, and every block $B$ from this family appears in any 
$x\in X$ with gaps ranging between $\tfrac{1-\frac 1r}{\mu(B)}$ and $\tfrac{1+\frac 1r}{\mu(B)}$, 
revealing strong repelling.
\endexample

\remark{Remark 1}
Viewing the blocks of length $nrN_0^{n-1}$ starting with a block from $\Cal C_1$ as a new alphabet, and 
repeating the above construction inductively, we can produce an example (with the measure of maximal entropy
on the intersection of systems created in consecutive steps) with entropy $\log_2N - 2\delta$, in which the 
strong repelling will occur with probability $\frac 1{n_k}$ for infinitely many lengths $n_k$.
\endremark

\remark{Remark 2}
The process described in the above remark is (somewhat coincidently; it was not designed for that)
bilaterally deterministic: for every $m\in\na$ the sigma-field $\Cal P^{(-\infty,-m]\cup[m,\infty)}$ 
equals the full (product) sigma-field. Indeed, suppose we see all entries of a $\Cal P$-name of a point $x$ 
except on the interval $(-m,m)$. In a typical point, this interval is contained between a pair of successive 
markers $i$ for some level $k$ of the inductive construction. Then, by examining this name's entries far enough 
to the left and right we will see complete all but one (the one covering the coordinate zero) blocks from the 
family $\Cal B_i$ which constitute the block $C\in\Cal C_i$ covering the considered interval. Because every 
block from $\Cal B_i$ is used in $C$ exactly once, by elimination, we will be able to determine the missing 
block from $\Cal B_i$ and hence all symbols in $(-m,m)$.
\endremark

\comment
\medskip
The next construction shows that there exists a process isomorphic to a Bernoulli process with an almost 
sure limit law $\tilde F\equiv 0$ for the normalized hitting times (strong attracting), achieved along a 
subsequence of upper density 1. In particular, this answers in the negative a question of Zaqueu Coelho 
([C]), whether all processes isomorphic to Bernoulli processes have necessarily the exponential limit 
law for the hitting (and return) times. 

\example{Example 2} We will build a decreasing sequence of subshifts of finite type (SFT's). In each we will 
regard the measure of maximal entropy. Begin with the full shift $X_0$ on a finite alphabet. Select $r$
words $W_1, W_2, \dots, W_r$ of some length $l$ and create $r$ SFT's: $X_0^{(1)}, X_0^{(2)},\dots, X_0^{(r)}$, 
forbidding one of these words in each of them, respectively. Choose a length $n$ so large that in the majority 
of blocks of this length in $X_0^{(i)}$ all words $W_j$ except $W_i$ will appear at least once. Now choose
another length $m$, such that in the majority of blocks of this length every block $C$ of length between $n$ 
and $n^2$ will appear many times. Now define $X_1$ as the subshift whose each point is a concatenation of the 
form $\dots B_1B_2\dots B_rB_1\dots$, where $B_i$ is a block appearing in $X_0^{(i)}$ of length either $m$ or 
$m+1$. Obviously, a typical block $C$ of any length between $n$ and $n^2$ appearing in $X_1$ comes from some 
$X_0^{(i)}$, hence contains all $W_j$'s except $W_i$, therefore in a typical $x\in X_1$, $C$ will appear many 
times within each component $B_i$ representing $X_0^{(i)}$, and then it will be absent for a long time, until 
the next representative of $X_0^{(i)}$. So, every such block will reveal strong attracting. It is not hard 
to see that $X_1$ is a mixing SFT and its d-bar distance from the full shift is small whenever the length 
$l$ of the (few) forbidden words $W_i$ is large. We can now repeat the construction starting with $X_1$, and 
radically increasing all parameters. We can arrange that the d-bar distances are summable, so the 
limit system $X$ (the intersection of the $X_k$'s), more precisely its measure of maximal entropy, is also a 
d@-bar limit. Each mixing SFT is isomorphic to a Bernoulli process and this property passes via d-bar limits 
(see [O], [Sh]), hence $X$ is also isomorphic to a Bernoulli process. This system has the almost sure limit law 
$\tilde F\equiv 0$ for hitting times (or $F\equiv 1$ for the return times) achieved along a \sq\ containing infinitely
many intervals of the form $[n,n^2]$. Such \sq\ has upper density 1. 
\endexample

\remark{Remark 3}
It is also possible to construct a process $X_h$ as above with any preassigned entropy $h$. On the other hand, 
it is well known ([Si]), that every measure-preserving transformation with positive entropy $h$ possesses a 
Bernoulli factor of the same entropy. By the Ornstein Theorem ([O]) this factor is isomorphic to 
$X_h$. The generator of $X_h$ appears as a partition of the space on which the initial measure-preserving 
transformation is defined. This proves the universality of {\it law of series}: {\bf in every 
measure-preserving transformation there exists a partition generating the full entropy, which has the
``strong repelling properties''} (i.e., almost sure limit law $\tilde F\equiv 0$ along a \sq\ of lengths of 
upper density 1). 
\endremark

\medskip
Various zero entropy processes with persistent repelling or attracting are implicit in the existing literature.
Extreme repelling (with intensity converging to $e^{-1}$ as the length of blocks grows) occurs for example 
in odometers, or, more generally, in rank one systems ([C-K]). For completeness, we sketch two zero entropy
processes with features of positive entropy: repelling, and the unbiased behavior.

\example{Example 3} Take the product of the independent Bernoulli process on two symbols with an odometer 
(the latter modeled by an adequate process, for example a regular Toeplitz subshift; see [D] for details on Toeplitz flows). 
Call this product process $X_0$. The odometer factor provides, for each $k\in\na$, markers dividing each element 
into so-called {\it $k$-blocks} of equal lengths $p_k$. Each $p_k$ is a multiple of $p_{k-1}$ and each $k$-block 
is a concatenation of $(k\!-\!1)$-blocks. Now we create a new process $X_1$ by ``stuttering'': if 
$x_0\in X_0$ is a concatenation $\dots ABCD\dots$ of $1$-blocks, we create $x_1\in X_1$ as 
$\dots AABBCCDD\dots$, with the number of repetitions 
$q_1=2$. In $X_1$ the lengths of the $k$-blocks for $k>1$ have doubled. Repeating the stuttering for 
$2$-blocks of $X_1$ with a number of repetitions $q_2\ge 2$, we obtain a process $X_2$. And so on. 
Because in each step we reduce the entropy by at least half, the limit process has entropy 
zero. If the $q_k$'s grow sufficiently fast, we obtain, like in the previous example, 
a system with strong attracting for a set of lengths of upper density 1. 
Consider a modification of this example where $q_k = 2$ for each $k$ and each pair $AA$ (also $BB$, 
etc.) is substituted by $A\overline A$ ($B\overline B$, etc.), where $\overline A$ is the ``mirror'' of $A$, 
i.e., with the symbols 0 are replaced by 1 and {\it vice versa}. It is not very hard to compute, that such process 
(although has entropy zero), has the same limit law properties as the independent process: almost sure 
convergence along the full \sq\ $(n)$ to the unbiased (exponential) limit law. 
\endexample

\remark{Remark 4} It is not hard to construct zero entropy processes with persistent mixed behavior. 
For example, applying the ``stuttering technique'' to an odometer one obtains a process in which a typical 
block $B$ occurs in periodically repeated pairs: $BB........BB........BB........$, i.e., with the function
$G_B\approx \min \{1,\frac t2\}$ (which reveals attracting with intensity 
$\frac{\log 2-1}2$ at $t_1=\log2$ and repelling with intensity $e^{-2}$ at $t_2=2$). We skip the details.
\endremark
\endcomment

\heading Questions \endheading
\remark{Question 1} Is there a speed of the convergence to zero of the joint measure of the ``bad'' blocks
in Theorem 1? More precisely, does there exist a positive function $s(n,\epsilon,\#\Cal P)$ converging to 
zero as $n$ grows, such that if for some $\epsilon$ and infinitely many $n$'s, the joint measure of the 
$n$-blocks which repel with intensity $\epsilon$ exceeds $s(n,\epsilon,\#\Cal P)$, then the process has 
necessarily entropy zero? (By the Example 1, $\frac1n$ is not enough.)
\endremark
\remark{Question 2}
Can one strengthen the Theorem 3 as follows: 
$$
\limsup_{n\to\infty} \,\tilde F_{x,n}\le 1-e^{-t} \ \ \mu\text{-almost everywhere?}
$$ 
\endremark
\remark{Question 3}
In Lemma 3, can one obtain $\Cal P^r$ conditionally $\beta$-independent of jointly the past and 
{\it all} return times $R^{(k)}_B$ ($k\ge 1$) (for sufficiently large $n$, with $\mu$-tolerance 
$\beta$ for $B\in\Cal P^{-n}$)? In other words, can the $\beta$-independent process \procbr\ be
obtained $\beta$-independent of the factor-process generated by the partition into $B$ 
and its complement?
\endremark
\remark{Question 4} (suggested by J-P. Thouvenot) Find a purely combinatorial proof of Theorem 1,
by counting the quantity of very long strings (of length $m$) inside which a positive fraction (in measure) 
of all $n$-blocks repel with a fixed intensity. For sufficiently large $n$ this quantity should 
be eventually (as $m\to\infty$) smaller than $h^m$ for any preassigned positive $h$.
\endremark
\comment
\remark{Question 5} As we have mentioned, we only know about conditions which ensure that the limit law
for the return time is exponential. It would be interesting to find a (large) class of positive entropy processes 
for which the distributions of return times are essentially deviated from exponential for bounded away from zero 
in measure collections of arbitrarily long blocks, i.e., a class of processes with persistent attracting. 
Can one prove that persistent attracting is, in some reasonable sense, a ``typical'' property in positive 
entropy, or that for a fixed measure-preserving transformation with positive entropy, a ``typical'' generator 
(partition) leads to persistent attracting?
\endremark
\endcomment

\enddocument